\documentclass[11pt]{amsart2000}
\newtheorem{theorem}{Theorem}
\newtheorem{fact}{Fact}

\begin{document}
\setlength{\unitlength}{0.01in}
\linethickness{0.01in}
\begin{center}
\begin{picture}(474,66)(0,0) 
\multiput(0,66)(1,0){40}{\line(0,-1){24}}
\multiput(43,65)(1,-1){24}{\line(0,-1){40}}
\multiput(1,39)(1,-1){40}{\line(1,0){24}}
\multiput(70,2)(1,1){24}{\line(0,1){40}}
\multiput(72,0)(1,1){24}{\line(1,0){40}}
\multiput(97,66)(1,0){40}{\line(0,-1){40}} 
\put(143,66){\makebox(0,0)[tl]{\footnotesize Proceedings of the Ninth Prague Topological Symposium}}
\put(143,50){\makebox(0,0)[tl]{\footnotesize Contributed papers from the symposium held in}}
\put(143,34){\makebox(0,0)[tl]{\footnotesize Prague, Czech Republic, August 19--25, 2001}}
\end{picture}
\end{center}
\vspace{0.25in}
\setcounter{page}{353}
\title{Nonstandard proofs of Eggleston like theorems}
\author{Szymon \.{Z}eberski}
\address{Institute of Mathematics, University of Wroc\l aw,
pl.\ Grunwaldzki 2/4, 50-384 Wroc\l aw, Poland}
\email{szebe@math.uni.wroc.pl}
\thanks{Szymon \.Zeberski,
{\em Nonstandard proofs of Eggleston like theorems},
Proceedings of the Ninth Prague Topological Symposium, (Prague, 2001),
pp.~353--357, Topology Atlas, Toronto, 2002}
\begin{abstract}
We prove theorems of the following form: if $A\subseteq {\mathbb R}^2$ is 
a ``big set'', then there exists a ``big set'' $P\subseteq {\mathbb R}$ 
and a perfect set $Q\subseteq {\mathbb R}$ such that 
$P\times Q\subseteq A$. 
We discuss cases where ``big set'' means: set of positive Lebesgue 
measure, set of full Lebesgue measure, Baire measurable set of second
Baire category and comeagre set. 
In the first case (set of positive measure) we obtain the theorem due to
Eggleston. 
In fact we give a simplified version of the proof given by J. Cicho\'{n}
in \cite{C}. 
To prove these theorems we use Shoenfield's theorem about absoluteness for 
$\Sigma^1_2$-sentences.
\end{abstract}
\maketitle

\section{Introduction}

We use standard set theoretical notation. 
By $\omega$ we denote the set of natural numbers. 
The cardinality of a set $X$ we denote by $|X|$.

By $[0,1]$ we denote the unit interval of the real line. 
By $\mbox{\sf Perf}([0,1])$ we denote the Polish space of all non-empty 
perfect subsets of the interval $[0,1]$ with the Hausdorff metric.
Similarly, by $\mbox{\sf Perf}({\mathbb R})$ we denote the Polish space
of all non-empty perfect subsets of the real line.
By $\mbox{\bf F}_\sigma $ we denote the family of all $F_\sigma $ subsets 
of the unit interval. 
By $\mbox{\bf G}_\delta $ we denote the family of all $G_\delta $ subsets
of the unit interval.

By $\lambda$ we denote Lebesgue measure on the real line, by $\lambda^2$ 
--- Lebesgue measure on the plane. 
${\mathbb L}$ denotes the $\sigma$-ideal of null sets 
(${\mathbb L} = \{A : \lambda (A)=0\}$). 
By ${\mathbb L}^+$ we denote the family of positive Lebesgue measure sets. 
${\mathbb K}$ denotes the $\sigma$-ideal of first Baire category sets.
${\mathbb K}^+$ denotes the family of Baire measurable sets of second
Baire category.
Let us recall that if $\mathcal{J}$ is a proper $\sigma$-ideal then 
$\mathcal{B} \subseteq \mathcal{J}$ is a base of $\mathcal{J}$ iff
$(\forall J\in \mathcal{J})(\exists B\in {\mathcal B})\ J\subseteq B$.
Furthermore, we consider a cardinal number 
$\mbox{\rm cof}(\mathcal{J})$ defined
as follows:
$$\mbox{\rm cof}(\mathcal{J}) = 
\min\{|\mathcal{B}| : \mathcal{B} \mbox{ is a base of } {\mathcal J}\}.$$

It is well known that we can extend every countable standard transitive 
model of {\rm ZFC} via forcing extension to a model of 
$\mbox{\rm ZFC} + 
(\mbox{\rm cof}({\mathbb L}) = \omega_1) + 
(2^\omega =\omega_2)$.
We can also extend every countable standard transitive model of
{\rm ZFC} via forcing extension to a model of
$\mbox{\rm ZFC} + 
(\mbox{\rm cof}({\mathbb K}) = \omega_1) +
(2^\omega =\omega_2)$.

We can code Borel subsets of the real line by functions from 
$\omega^\omega$ in an absolute way (see \cite{J}). 
Functions which code Borel sets are called Borel codes.
We use the folowing notation: 
if $f\in\omega^\omega $ is a Borel code, then
$\# f\subseteq {\mathbb R}$ is a set coded by $f$.

By a canonical Polish space, we understand a countable product of spaces
$\{0,1\}^{\omega }$, $\omega^\omega $, $[0,1]$, $\mbox{\sf Perf}([0,1])$ 
and so on. 
A sentence $\varphi$ is a $\Sigma_{2}^{1}$-sentence if for some canonical
Polish spaces $X$, $Y$ and some Borel subset $B\subseteq X\times Y$ we 
have
$\varphi =(\exists x \in X)(\forall y \in Y)\ (x,y) \in B$. 
Spaces $X$, $Y$ and the set $B$ are called ``parameters'' of the sentence
$\varphi$. 
We shall use the following classical theorem about absoluteness of
$\Sigma_{2}^{1}$-sentences (see \cite{SH}):

\begin{theorem}[Shoenfield]
Suppose that $M\subseteq N$ are standard transitive models of the theory
ZF such that $\omega _{1}^{N}\subseteq M$. 
Let $\varphi$ be a $\Sigma_{2}^{1}$-sentence with parameters from the
model $M$. 
Then
$$M\models\varphi\ \iff\ N\models\varphi.$$
\end{theorem}

Notice that if the model $N$ is a generic extension of the standard 
transitive model $M$ then both models $M$ and $N$ have the same ordinal 
numbers, so the inclusion $\omega _{1}^{N}\subseteq M$ trivially holds.

Let us note that if $A$ is a coanalytic subset of a Polish space and
$|A|>\omega_1$ then there exists a nonempty perfect subset of $A$. 
Indeed, there exists a family $\{B_\alpha \}_{\alpha <\omega_1}$ of Borel
sets such that $A=\bigcup_{\alpha <\omega_1}B_\alpha$ (see \cite{CKW}). 
Hence, if $|A|>\omega_1$ then there exists $\alpha <\omega_1$ such that 
$|B_\alpha |>\omega_1$, so $B_\alpha $ has a nonempty perfect subset. 

\section{Main result}

We give a new proof of Eggleston's theorem (see \cite{E}) and then by a
simple modification of this proof we get new results about inscribing
special rectangles into ``big sets'' in the sense of measure or category.

\subsection{Measure case}

In the paper \cite{CKP} we can find the following result:

\begin{theorem}[Cicho\'{n}, Kamburelis, Pawlikowski]
There exists a dense subset of the measure algebra $\mbox{\sf Borel}
({\mathbb R})/{\mathbb L}$ of cardinality $\mbox{\rm cof}({\mathbb L})$.
\end{theorem}

Using this theorem we can deduce the following fact:

\begin{fact}\label{b}
There exists a family 
${\mathcal B} \subseteq \mbox{\sf Perf}({\mathbb R})\cap {\mathbb L}^+$
such that $|{\mathcal B}|=\mbox{\rm cof}({\mathbb L})$
and
$$(\forall A\subseteq {\mathbb R})\
\big( 
\lambda (A)>0
\longrightarrow
(\exists P\in {\mathcal B})\ P\subseteq A
\big)
.$$
\end{fact}

Using this fact, we can obtain a simple proof of the following theorem from
\cite{E}:

\begin{theorem}[Eggleston]\label{e}
Let $A$ be a subset of the plane of positive Lebesgue measure. 
Then we can find perfect subsets of the real line $P$, $Q$ such that
$\lambda (P)>0$ and $P\times Q\subseteq A$.
\end{theorem}

\begin{proof}
Without loss of generality we can assume that $A$ is a Borel set.
Let $V'$ be a generic extension of the universe $V$ such that 
$$V'\models \mbox{\rm ZFC} + 
(\mbox{\rm cof}({\mathbb L})=\omega_1) +
(2^\omega =\omega_2).$$
For a while we will work in the universe $V'$. 
For $y\in {\mathbb R}$ let 
$A^y=\{x\in {\mathbb R} : (x,y)\in A\}$
and we put 
$$Y=\{y\in {\mathbb R} : \lambda (A^y)>0\}.$$
It follows from Fubini's theorem that $\lambda (Y)>0$.
In particular $|Y|=\omega_2$.
Using Fact \ref{b} we get a family
${\mathcal B}\subseteq \mbox{\sf Perf}({\mathbb R})\cap {\mathbb L}^+$
such that
$$
(\forall A\subseteq {\mathbb R})\ 
\big( 
\lambda (A)>0
\longrightarrow
(\exists P\in {\mathcal B})\ P\subseteq A
\big)
$$
and $|{\mathcal B}|=\mbox{\rm cof}({\mathbb L})=\omega_1$.
For every $y\in Y$ we can find $P\in {\mathcal B}$ such that 
$P\subseteq A^y$.
So there exists $P\in\mbox{\sf Perf}({\mathbb R})\cap {\mathbb L}^+$ such
that 
$|\{y\in Y : P\subseteq A^y\}|=\omega_2$.
But $\{y\in Y : P\subseteq A^y\}$ is a coanalytic set of cardinality 
$\omega_2$ so we can find a perfect set $Q$ such that
$Q\subseteq\{y\in Y : P\subseteq A^y\}$.
It means that in the universe $V'$ the following sentence holds:
$$
(\exists P,Q\in \mbox{\sf Perf}({\mathbb R}))
(\forall r\in {\mathbb R}^2)\
\lambda (P)>0 \wedge 
(r\in P\times Q\rightarrow r\in A)
.
$$
But this is a $\Sigma_2^1$-sentence, hence it is also true in the original 
universe V.
\end{proof}

Now, we show the theorem about inscribing special rectangles into sets of
full Lebesgue measure. At first, we have to formulate the following fact:

\begin{fact}\label{f}
There exists a family ${\mathcal F}\subseteq \mbox{\bf F}_\sigma$ of
measure one subsets of the unit interval such that
$|{\mathcal F}|=\mbox{\rm cof}({\mathbb L})$ and
$$
(\forall A\subseteq [0,1])\
\big(
\lambda (A)=1
\longrightarrow
(\exists F\in {\mathcal F})\ F\subseteq A
\big).
$$
\end{fact}

Using this fact, we prove the following theorem:

\begin{theorem}\label{mf}
Let $A$ be a subset of $[0,1]^2$ such that $\lambda^2 (A)=1$. 
Then we can find two sets $F$, $Q\subseteq [0,1]$ such that $F$ is a
$F_\sigma $ set, $\lambda (F)=1$, $Q$ is a perfect set and 
$F\times Q\subseteq A$.
\end{theorem}

\begin{proof}
Without loss of generality $A$ is a Borel set.
Let $V'$ be a generic extension of the universe $V$ such that 
$$
V'\models \mbox{\rm ZFC} + 
(\mbox{\rm cof}({\mathbb L}) = \omega_1) +
(2^\omega = \omega_2).
$$
Similarly as in the proof of theorem \ref{e}, using Fact \ref{f}, we 
deduce that in the universe $V'$ the folowing sentence holds
$$
(\exists F\in \mbox{\bf F}_\sigma)(\exists Q\in \mbox{\sf Perf}([0,1]))\
\lambda (F)=1 
\wedge 
F\times Q\subseteq A
.
$$
But we can formulate this sentence in the folowing way
$$
(\exists f\in \omega^\omega)(\exists Q\in \mbox{\sf Perf}([0,1]))\
f \mbox{ codes $F_\sigma$ set }
\wedge 
\lambda (\# f)=1
\wedge 
\# f\times Q\subseteq A.
$$
This is also a $\Sigma_2^1$-sentence, hence it is true in the original
universe $V$.
\end{proof}

\subsection{Category case}

Now, we will formulate and prove theorems about inscribing special 
rectangles into ``big sets'' in the sense of category.
We start with two observations (similar to Fact \ref{b} and \ref{f}).
Namely, we have the following facts:

\begin{fact}\label{g}
There exists a family 
${\mathcal G}\subseteq \mbox{\bf G}_\delta \cap {\mathbb K}^+$
such that $|{\mathcal G}|=\mbox{\rm cof}({\mathbb K})$ and
$$(\forall A\in {\mathbb K}^+)(\exists G\in {\mathcal G})\ G\subseteq A.$$
\end{fact}

\begin{fact}\label{gc}
There exists a family ${\mathcal C}\subseteq \mbox{\bf G}_\delta $ of
comeagre subsets of the unit interval such that 
$|{\mathcal C}|=\mbox{\rm cof}({\mathbb K})$ and
$$
(\forall A\subseteq [0,1])\
\big(
A \mbox{ is comeagre}
\longrightarrow
(\exists C\in {\mathcal C})\ C\subseteq A
\big)
.
$$
\end{fact}

Using these facts, we prove the following theorems:

\begin{theorem}\label{gt}
Let $A\subseteq [0,1]^2$ be a Baire measurable set of second Baire category.
Then we can find two sets $G$, $Q\subseteq [0,1]$ such that $G$ is a
$G_\delta$ set of second Baire category, $Q$ is a perfect set and 
$G\times Q\subseteq A$.
\end{theorem}

\begin{theorem}\label{gct}
Let $A\subseteq [0,1]^2$ be a comeagre set. 
Then we can find two sets $C$, $Q\subseteq [0,1]$ such that $C$ is a
comeagre $G_\delta$ set, $Q$ is a perfect set and $C\times Q\subseteq A$.
\end{theorem}

\begin{proof}[Proof of Theorem \ref{gt} and Theorem \ref{gct}]
We start with extending the universe $V$ to $V'$ (by generic 
extension) such that 
$$V'\models \mbox{\rm ZFC} + 
(\mbox{\rm cof}({\mathbb K}) = \omega_1) +
(2^\omega =\omega_2).$$
Similarly as in the proof of Theorem \ref{e} (naturally using Facts 
\ref{g}, \ref{gc} respectively for Theorem \ref{gt} and \ref{gct}), we
deduce that if $A$ is a Baire measurable second Baire category set 
(Theorem \ref{gt}) then in the model $V'$ the following sentence holds:
$$
(\exists G\in \mbox{\bf G}_\delta)(\exists Q\in \mbox{\sf Perf}([0,1]))\
G\in {\mathbb K}^+ 
\wedge 
G\times Q\subseteq A
.
$$
When $A$ is comeagre (Theorem \ref{gct}) then in the universe $V'$ the 
following sentence holds:
$$
(\exists C\in \mbox{\bf G}_\delta)(\exists Q\in \mbox{\sf Perf}([0,1]))\
C \mbox{ is comeagre } 
\wedge 
C\times Q\subseteq A
.
$$
Similarly as in the proof of Theorem \ref{mf}, we can use Borel coding
for replacing quantifiers of the type $(\exists A\in \mbox{\bf G}_\delta)$
by 
$$
(\exists a\in \omega^\omega) \ a \mbox{ codes } G_\delta\mbox{ set}
$$ 
and replace $A$ by $\# a$. 
After this modification we obtain equivalent sentences, which are
$\Sigma_2^1$-sentences, hence they are true in the original universe $V$.
\end{proof}

\providecommand{\bysame}{\leavevmode\hbox to3em{\hrulefill}\thinspace}
\providecommand{\MR}{\relax\ifhmode\unskip\space\fi MR }
\providecommand{\MRhref}[2]{%
  \href{http://www.ams.org/mathscinet-getitem?mr=#1}{#2}
}
\providecommand{\href}[2]{#2}

\end{document}